\documentstyle[12pt]{article}
\begin{document}
\renewcommand{\thefootnote}{\fnsymbol{footnote}}

\begin{center}
{\normalsize\bf
$\sigma$ - models on the quantum group manifolds $SL_{q}(2,R)$,
$SL_{q}(2,R)/U_{h}(1)$, $C_{q}(2|0)$ and infinitesimal
transformations.\\
} \vspace{0.20cm} V.D.  Gershun \vspace{0.20cm}\\ {\it Kharkov Institute
of Physics and Technology, 310108, Kharkov, Ukraine}\\ \vspace{0.20cm}
\end{center}

\begin{quotation}
{\small\rm
 The differential and variational calculus on the $SL_{q}(2,R)$ group
is constructed. The spontaneous breaking symmetry in the WZNW model with
$SL_{q}(2,R)$ quantum group symmetry and in the $\sigma$-models with
${SL_{q}(2,R)/U_{h}(1)}$ ,$C_{q}(2|0)$ quantum group symmetry is
considered. The Lagrangian formalism over the quantum group manifolds is
discussed. The classical solution of $C_{q}(2|0)$ {$\sigma$}-model is
obtained.}
\end{quotation}
\renewcommand{\thefootnote}{\arabic{footnote}}
\setcounter{footnote}0

\vspace{0.20cm}

{\bf 1. Differential calculus on the $SL_{q}(2,R)$ group.}\\ The matrix
quantum group [1] $G= SL_{q}(2,R)$ is defined by the q-commutation
relations (C.R.) of its group parameters. Let $$g=\left(
\begin{array}{cc} a^1&a^2\\ a^3&a^4 \end{array} \right),\;
\begin{array}{ccc} a^1a^2 = qa^2a^1,& a^2a^4 = qa^4a^2,& a^2a^3 = a^3a^4
\\ a^1a^3 = qa^3a^1,& a^3a^4 = qa^4a^3,& a^1a^4 = a^4a^1 + ( q - q^{-1})
a^2a^3 \end{array} \eqno(1) $$\\ ${a^k}$ - hermitian, $|q|=1, Det_{q}g=
a^1a^4- qa^2a^3= 1$.\\  For any elements $g$, $g'$ $\in$ $SL_{q}(2,R)$
element $g'' = {g'}{g}$ will belong $SL{q}(2,R)$ if
$a^{k '}a^{l} = a^{l}a^{k '}$. In the Gauss decomposition [2]
 $$g= \left( \begin{array}{cc} 1&\varphi_{-}\\ 0&1 \end{array} \right)
 \left( \begin{array}{cc} 1&0\\ \varphi_{+}&1 \end{array} \right) \left(
\begin{array}{cc}
\rho&0\\
0&\rho^{-1}
\end{array}
\right)=
\left(
\begin{array}{cc}
{\rho + \varphi_{-}\varphi_{+}\rho}&\varphi_{-}\rho^{-1}\\
\varphi_{+}\rho&\rho^{-1}
\end{array}
\right)
\eqno(2)
$$
the C.R. are:
$$\rho\varphi_{\pm}=q\varphi_{\pm}\rho,\;\;\varphi_{-}\varphi_{+}=
q^{2}\varphi_{+}\varphi_{-}
\eqno(3)
$$
  Let the quantum group is a manifolds of any possible transformations
$g' = g g_{0}$. There are two kinds of the variation: the variation
in the neighborhood of the arbitrary point of the group space $g'= g
+ dg$ and variation in the neighborhood of the unit of the group $g = 1
+
\delta g$. First variation defines the group invariants: element of the
distance between two neighboring points, element of the volume around
the
point. Second variation defines the group symmetry of this invariants.
  The C.R. between variation $dg$ and $g$ define the type of the
differential calculus.
      The left-invariant differential calculus [3] on the
$GL_{q}(2,C)$ group, matched with the differential calculi on the
$SL_{q}(2,C)$ subgroup and on the Borel subgroups $B_{L}(C), B_{U}(C)$,
was constructed in [4,5].  Let $\omega=g^{-1}dg$ is the
left differential Kartan 1-form  $$ \omega= \left(
\begin{array}{cc} \omega^{1}&\omega^{2}\\ \omega^{3}&\omega^{4}
\end{array} \right),\;\; Tr_{q}\omega=q^{2}\omega^{1}+\omega^{4}=0
\eqno(4)
$$
 The differential calculus on the $SL_{q}(2,R)$ group is defined by the 
C.R.
 $$ \begin{array}{ccc}
\omega^{1}\rho= {1\over q^{2}}\rho\omega^{1},&
\omega^{2}\rho= {1\over q}\rho\omega^{2},&
\omega^{3}\rho= {1\over q}\rho\omega^{3}\\
\omega^{1}\varphi_{\pm}= \varphi_{\pm}\omega^{1},&
\omega^{2}\varphi_{\pm}= \varphi_{\pm}\omega^{2},&
\omega^{3}\varphi_{\pm}=\varphi_{\pm}\omega^{3}
\end{array}
\eqno(5)
$$
 The C.R. between the group parameters and their differentials are more
 complicated:
$$
\begin{array}{cc}
{d\rho \rho={1\over q^{2}}\rho d\rho},&
{d\varphi_{+} \varphi_{+}={1\over q^{2}}\varphi_{+} d\varphi_{+}+
(q^{4}-1)\varphi_{+}^{3} d\varphi_{-}}\\
{d\varphi_{-} \varphi_{-}=q^{2}\varphi_{-} d\varphi_{-}},&
{d\rho \varphi_{-}=q\varphi_{-} d\rho}
\end{array}
$$
$$
\begin{array}{ccc}
d\varphi_{-}\varphi_{+}=q^{2}\varphi_{+}d\varphi_{-},&
d\varphi_{+}\varphi_{-}={1\over q^{2}}\varphi_{-}d\varphi_{+},&
d\rho d\varphi_{-}=-qd\varphi_{-}d\rho\\
d\varphi_{-}d\varphi_{+}=-q^{2}d\varphi_{+}d\varphi_{-},&
d\varphi_{-}\rho={1\over q}\rho d\varphi_{-}
\end{array}
\eqno(6)
$$
$$
\begin{array}{cc}
d\varphi_{+}\rho={1\over q}\rho d\varphi_{+}-
q(q^{2}-1)\varphi_{+}^{2}\rho d\varphi_{-},\;\;\;
d\rho\varphi_{+}=q\varphi_{+}d\rho - q^{2}(q^{2}-1)\varphi_{+}^{2}
\rho d\varphi_{-}\\
d\rho d\varphi_{+}+qd\varphi_{+}d\rho+q^{3}(q^{2}-1)
\varphi_{+}^{2}d\varphi_{-}d\rho-{(q^{4}-1)\over q^{3}}
\varphi_{+}\rho d\varphi_{-}d\varphi_{+}=0
\end{array}
$$
 The Kartan 1-forms are:
$$
{\omega^{1}=\rho^{-1}d\rho+\varphi_{+}d\varphi_{-}},\;\;
{\omega^{3}={1\over q}\rho^{2}d\varphi_{+}-q^{5}\varphi_{+}^{2}
\rho^{2}d\varphi_{-}},\;\;
{\omega^{2}=q\rho^{-2}d\varphi_{-}}
\eqno(7)
$$
$$
\begin{array}{cc}
(\omega^1)^2 = (\omega^2)^2 = (\omega^3)^2 = 0 &\;\;\;
\omega^4 = -q^2\omega^1\\ \omega^1\omega^2 + q^4\omega^2\omega^1 =
0&\;\;\; \omega^1\omega^3 + q^{-4}\omega^3\omega^1 = 0\\
\omega^2\omega^3
+ q^{-2}\omega^3\omega^2 = 0 \end{array} $$ The left vector fields
  $\nabla_{k}$ can be obtained from the applying the left differential
to
an arbitrary function on the quantum group $df =
(f{\partial\over\partial
a^k})da^k = (f\nabla_{k})\omega^{k}$.  $$ \begin{array}{cc} \nabla =
\left(\begin{array}{cc}
\nabla_{1}&\nabla_{2}\\
\nabla_{3}&\nabla_{4}
\end{array}
\right)&\;\;\hat{\nabla_{1}} = \nabla_{1}-q^2\nabla_{4}\;\;\;\;\;\;
\hat{\nabla_{4}} = \nabla_{1}+\nabla_{4}
\end{array}
\eqno(8)
$$
  The C.R. for vector fields have following form
$$
\begin{array}{ccc}
\rho\hat{\nabla_{1}} = {1\over q^2}\hat{\nabla_{1}}\rho + \rho &\;\;
\varphi_{-}\hat{\nabla_{1}} = \hat{\nabla_{1}}\varphi_{-}&\;\;
\varphi_{+}\hat{\nabla_{1}} = \hat{\nabla_{1}}\varphi_{+}\\
\rho\nabla_{2} = {1\over q}\nabla_{2}\rho - \varphi_{+}{\rho}^3 &\;\;
\varphi_{-}\nabla_{2} = \nabla_{2}\varphi_{-} + {1\over q}\rho^2 &\;\;
\varphi_{+}\nabla_{2} = \nabla_{2}\varphi_{+} + q\varphi_{+}^2\rho^2\\
\rho\nabla_{3} = {1\over q}\nabla_{3}\rho &\;\;
\varphi_{-}\nabla_{3} = \nabla_{3}\varphi_{-} &\;\;
\varphi_{+}\nabla_{3} = \nabla_{3}\varphi_{+} + q\rho^{-2}\\
\end{array}
\eqno(9)
$$
$$
q^2\hat{\nabla_{1}}\nabla_{3}-{1\over q^2}\nabla_{3}\hat{\nabla_{1}}=
(q^2+1)\nabla_{3},\;q^2\nabla_{2}\hat{\nabla_{1}}-{1\over
q^2}\hat{\nabla_{1}}\nabla_{2}=(q^2+1)\nabla_{2},\;
\nabla_{3}\nabla_{2}-{1\over
q^2}\nabla_{2}\nabla_{3}=\hat{\nabla_{1}}
\eqno(10)
$$
and $\nabla_{k}$ have the form
$$
\hat{\nabla_{1}} = {\partial\over \partial\rho}\rho,\;\;\nabla_{2} =
{1\over q}{\partial\over \partial\varphi_{-}}\rho^2 - {\partial\over
\partial\rho}\varphi_{+}\rho^3 + q{\partial\over
\partial\varphi_{+}}\varphi_{+}^2\rho^2,\;\;\nabla_{3} = q{\partial\over
\partial\varphi_{+}}\rho^{-2}
$$
  The left vector fields and the left derivatives act on the any
function
of the group parameters from the right side.\\
{\bf 2.WZNW model on the $SL_{q}(2,R)$ group.}\\ The existing of
the quantum group structure in the WZNW model was shown in [6,7]. The
$\sigma$-models with a quantum group symmetry was
considered in [2,8,9,10]. To construct the WZNW model with $SL_{q}(2,R)$
group symmetry, we consider the space ${M^{1,1}\oplus SL_{q}(2,R)}$,
where
$M^{1,1}$ is the commutative (undeformed) space.  The element of the
volume in $M^{1,1}$ space,which is the invariant of $SL_{q}(2,R)$, is
$${Tr_{q}[\omega(d)\wedge dz^{\mu}][\omega(d)\wedge dz^{\mu}]\over
{2\epsilon ^{\lambda\rho} dz^{\lambda}\wedge dz^{\rho}}}=
Tr_{q}(\omega_{\mu}\omega^{\mu}) d^{2}z , \eqno(11) $$ where
$\omega(d)=\omega_{\mu} dz^{\mu}, z^{\mu}\epsilon M^{1,1},
\mu=1,2.$ For any $2\times 2$ matrix $A, Tr_{q}A=q^{2}A^{1}+A^{4}$.As
a result we have
$$ Tr_{q}(\omega_{\mu}\omega^{\mu})=q^{5}[2]_{q}
\rho^{-2}\partial_{\mu}\rho\partial^{\mu}\rho+q^{5}
[2]_{q}\rho^{-1}\varphi_{+}(\partial_{\mu}\varphi_{-}
\partial^{\mu}\rho+{1\over q}
\partial_{\mu}\rho\partial^{\mu}\varphi_{-})+
$$
$$
(\partial_{\mu}\varphi_{-}\partial^{\mu}\varphi_{+}+
q^{2}\partial_{\mu}\varphi_{+}\partial^{\mu}\varphi_{-})-
q^{2}(q^{4}-1)\varphi_{+}^{2}\partial_{\mu}\varphi_{-}
\partial^{\mu}\varphi_{-}
\eqno(12)
$$
 The C.R. are now in the same space-time point ,
$d\rho=\partial_{\mu}\rho dz^{\mu}, d\varphi_{\pm}=
\partial_{\mu}\varphi_{\pm} dz^{\mu}$ and
$[n]_{q}={{q^{n}-q^{-n}}\over{q-q^{-1}}}$. The Wess-Zumino term
$$
Tr_{q}(\omega (d)\wedge\omega (d)\wedge \omega (d))=
{q[2]_{q}[3]_{q}\over 6}\varepsilon^{\mu\nu\lambda}
\partial_{\lambda}(\rho^{-1}\partial_{\mu}\rho\partial_{\nu}\varphi_{-}
\varphi_{+}) d^{3}z
\eqno(13)
$$
is the total derivative.
Finally, the WZNW-action  with the $SL_{q}(2,R)$ quantum group symmetry
describes the 2-dimensional relativistic string in the background
gravity
and antisymmetric fields
$$ S[\rho,\varphi_{-},\varphi_{+}]={k\over
4\pi}\int d^{2}z(G_{AB}\partial_{\mu}X^{A}\partial^{\mu}X^{B}+
B_{AB}\varepsilon_{\mu\nu}\partial^{\mu}X^{A}\partial^{\nu}X^{B}) ,
\eqno(14)
$$
where $X^{A}=(\rho,\varphi_{-},\varphi_{+})$ and the background gravity
and
antisymmetric fields have the following form:
$$
G_{AB}=
\left(
\begin{array}{ccc}
q^{5}[2]_{q}\rho^{-2}&q^{4}[2]_{q}\rho^{-1}\varphi_{+}&0\\
q^{5}[2]_{q}\rho^{-1}\varphi_{+}&-q^{2}(q^{4}-1)\varphi_{+}^{2}&1\\
0&q^{2}&0
\end{array}
\right),
B_{AB}=
{q^{3}[2]_{q}[3]_{q}\over 6}\varphi_{+}\rho^{-1}
\left(
\begin{array}{ccc}
0&1&0\\
-1&0&0\\
0&0&0
\end{array}
\right)
$$
  The group symmetry of this model is $SL_{q}(2,R)\otimes SL_{q}(2,R)$,
because under the left multiplication on the group
$g'=g_{0}g$ the differential forms of Kartan are invariant, $\omega
^\prime=\omega$, and under the right multiplication $g'=gg_{0}$
the differential forms are covariant, $\omega ^{\prime}=
{g_{0}^{-1}}\omega
g_{0}$.But $Tr_{q}A$ is invariant of the transformation
$A'=g_{0}^{-1}Ag_{0}$, because the elements of matrix $A$ commute with
the
elements of matrix $g_{0}$, by definition of the quantum
group.Therefore,
this model describes the spontaneous breaking of the $SL_{q}(2,R)\otimes
SL_{q}(2,R)$ symmetry to the $SL_{q} (2,R)$ one.\\ {\bf 3.$\sigma$-model
on the $SL_{q}(2,R)/U_{h}(1)$ group.}\\ Let us consider the spontaneous
breaking symmetry in the $\sigma-$ model with the \\
$SL_{q}(2,R)/U_{h}(1)$
group symmetry.  Let $G=KH$, K-coset, H-subgroup.\\ The Kartan 1-forms
$$
k^{-1}dk= \left( \begin{array}{cc}
q^{2}\varphi_{+}d\varphi_{-}&d\varphi{-}\\
d\varphi_{+}-q^{2}\varphi_{+}^{2}d\varphi_{-}&
-\varphi_{+}d\varphi_{-}
\end{array}
\right)=\omega +\theta ,
\eqno(15)
$$
where $\omega\epsilon K, \theta\epsilon H$ and the coset elements
$\varphi_{\pm}$ commute with the subgroup parameter $\rho$ and satisfy
to
C.R. of $SL_{q}(2,R)$ group among themselves.  There is a question:how
do
coset and subgroup separate from $k^{-1}dk$?  In opposite to the
classical case, there is the 3-parametric family of the $U(1)$
subgroups. The Lagrangian has
the following form:
$$
L_{n}={1\over 2}Tr_{q}(\omega_{\mu}\omega^{\mu})= {(q^{4}+1)\over
4q^{4}}(\partial_{\mu}\varphi_{-}\partial^{\mu}
\varphi_{+}+q^{2}\partial_{\mu}\varphi_{+}\partial^{\mu}\varphi_{-})
-c_{n}(q)\varphi_{+}^{2}\partial_{\mu}\varphi_{-}\partial
^{\mu}\varphi_{-},
\eqno(16)
$$ where $c_{n}(q)$ depends on the choice of a subgroup. There are three
most interesting examples. \\ {\bf 1)\/}Undeformed $U(1)$ subgroup:
$c_{1}={{2q^{4}-q^{2}+1}\over 2}$ $$ \omega=\left( \begin{array}{cc}
 (q^{2}-1)\varphi_{+}d\varphi_{-}&d\varphi_{-}\\
d\varphi_{+}-q^{2}\varphi_{+}^{2}d\varphi_{-}&0
\end{array}
\right),\;\;\;\theta=\varphi_{+}d\varphi_{-}
\left(
\begin{array}{cc}
1&0\\
0&-1
\end{array}
\right)
\eqno(17)
$$
  The algebra symmetry of this Lagrangian is defined by the
Maurer-Kartan
equations:
$$
d\theta=-
\left(
\begin{array}{cc}
q^{-2}&0\\
0&1
\end{array}
\right)
\omega\omega+(q^{2}-1)\omega\theta , d\omega=- \left(
\begin{array}{cc} {q^{2}-1}\over{q^{2}}&0\\ 0&0 \end{array} \right)
\omega\omega-q^{3}[2]_{q}\omega\theta , \theta\omega=q^{4}\omega\theta
$$
The C.R. between
the coset and the subgroup forms are common for all of the examples
$$
\begin{array}{cc} {\omega ^{1}\omega ^{3}+q^{4}\omega ^{3}\omega ^{1}=0}
,
& {\omega ^{2}\omega ^{3}+q^{2}\omega ^{3}\omega ^{2}=0}\\ {\omega
^{1}\omega ^{3}+q^{4}\omega ^{3}\omega ^{1}=0} , & {\omega ^{4}\omega
^{3}+q^{4}\omega ^{3}\omega ^{4}=0} \end{array}
\eqno(18) $$\\ {\bf
2)\/}Classical coset structure: $c_{2}={{q^{6}+1}\over 4}$ $$
\omega=\left( \begin{array}{cc} 0&d\varphi_{-}\\
d\varphi_{+}-q^{2}\varphi_{+}^{2}d\varphi_{-}&0
\end{array}
\right),\;\;\;\theta=\varphi_{+}d\varphi_{-}
\left(\begin{array}{cc}
q^{2}&0\\
0&-1
\end{array}
\right)
\eqno(19)
$$
$$
d\theta=-\omega\omega,\;\;
d\omega=-\omega\theta-\theta\omega ,\;\theta\omega =
q^{2}\omega\theta $$ {\bf 3)\/}There is one of the examples of the 2-
parametric family $U_{q}(1)$ subgroups:
$$
\begin{array}{ccc}c_{3}={{2q^{4}-q^{2}+1}\over 2q^{2}},& \omega=\left(
\begin{array}{cc} {(q^{2}-1)\over
q^{2}}\varphi_{+}d\varphi_{-}&d\varphi_{-}\\
d\varphi_{+}-q^{2}\varphi_{+}^{2}d\varphi_{-}&0
\end{array}
\right),&\theta={1\over q^{2}}\varphi_{+}d\varphi_{-}
\left(\begin{array}{cc}
1&0\\
0&-q^{2}
\end{array}
\right)
\end{array}
\eqno(20)
$$
$$
d\theta=-
\left(
\begin{array}{cc}
q^{-4}&0\\
0&1
\end{array}
\right)
\omega\omega+(q^{4}-1)\omega\theta,\;
d\omega=-
\left(
\begin{array}{cc}
{{q^{4}-1}\over{q^{4}}}&0\\
0&0
\end{array}
\right)
\omega\omega-q^{4}[2]_{q}\omega\theta,\theta\omega=q^{6}\omega\theta
$$
  Why we have obtained different algebras of a symmetry for the same
subgroup?
  That is possible because we can use the different map from the algebra
to
the group , for example:\\
$$g=\exp(\varphi_{-}\tau_{+})\exp(\varphi_{+}\tau_{-})\exp({\ln\rho}
\tau_{3}) ,
\eqno(21)
$$
where $\tau$ are the Pauli matrices -- the fundamental representation of
the $U_{q}(SL(2,R)$ algebra.The group stability of the vacuum is $U(1)$.
In the another parametrization\\
$$g=\exp(\varphi_{-}\tau_{+})\exp(\varphi_{+}\tau_{-})(1-{(q^{2}-1)\over
q^{2}}\nabla_{3})^{\ln\rho\over \ln q^{-2}},\; \nabla_{3}= \left(
\begin{array}{cc}
1&0\\
0&-q^{2}
\end{array}
\right)
\eqno(22)
$$
the group stability of the vacuum is $U_{q}(1)$.
 The group symmetry of this Lagrangians is $SL_{q}(2,R)$ spontaneously
breaked to $U_{h}(1),\;\;h=q^{\pm 2n},n=0,1... $.  Under the left
multiplication on the group $G'=G_{0}G$,the differential form
$G^{-1}dG=H^{-1}(\omega +\theta )H={G'} ^{-1}dG'$.  Therefore,
$\omega^{\prime} +\theta^{\prime} =H'H^{-1}(\omega +\theta )H{H'}^{-1}$
.
  Again, the decomposition on the coset and
the subgroup forms is not unique after transformation . The group
transformation can transform the Lagrangian with the $U_{h_{1}}(1)$
subgroup of the vacuum stability to the Lagrangian with the
$U_{h_{2}}(1)$
subgroup.\\
{\bf 4.Variational calculus on the $SL_{q}(2,R)$ group.}\\
  It is possible to obtain  the variational calculus on the group by two
ways: from the C.R. between the left vector fields and group parameters
and from the infinitesimal transformations on the group. Let us multiply
the C.R. (8) between $\nabla_{n}$ and group parameters on the parameters
of transformation $R^n$. The form of the infinitesimal transformations
of
the group parameters is obtained under the requirement
$$
[X_{A},\nabla_{n}R^n] =X_{A}\delta_{R^n},\;\;\;X_{A}=(\rho,
\varphi_{-}, \varphi_{+}),\;\;\; [A,B]=AB-BA \eqno(23) $$ By imposing
the
 C.R. between the parameters of infinitesimal transformations and group
parameters $$ \begin{array}{ccc} \rho R^{1} = q^2
R^{1}\rho&\;\;\varphi_{-}R^{1} = R^{1}\varphi_{-}&\;\;\varphi_{+}R^{1} =
R^{1}\varphi_{+}\\ \rho R^{2} = q R^{2}\rho&\;\;\varphi_{-}R^{2} =
R^{2}\varphi_{-}&\;\;\varphi_{+}R^{2} = R^{2}\varphi_{+}\\
\rho R^{3} = q R^{3}\rho&\;\;\varphi_{-}R^{3} =
R^{3}\varphi_{-}&\;\;\varphi_{+}R^{3} = R^{3}\varphi_{+}
\end{array}
\eqno(24)
$$
we obtain the infinitesimal transformation of the group parameters
$$
\begin{array}{ccc}\rho\delta = \rho R^{1} -
\varphi_{+}\rho^3R^{2}&\;\;\varphi_{-}\delta = {1\over q}\rho^2
R^{2}&\;\;\varphi_{+}\delta = q\varphi_{+}^2\rho^2 R^{2} +
q\rho^{-2}R^{3}
\end{array}
\eqno(25) $$ The same result we can obtain from the right infinitesimal
 multiplication on the group $g'= g g_{0}$, where $g_{0}=1 + \delta
g_{0}$.  For $$ \delta g_{0} = \left(\begin{array}{cc} R^{1}&R^{2}\\
R^{3}&-q^2R^{1}
\end{array}
\right)
\eqno(26)
$$
we see, that $dg = g\delta g_{0}$ and C.R. for $\delta g_{0}$ are the
same
as for left forms $\omega$ simultaneously with condition $R^{4} = -q^2
R^{1}$. The C.R. of the variational calculus
$$
\begin{array}{ccc}
(\rho\delta) \rho = {1\over q^2}\rho
(\rho\delta)&\;\;(\varphi_{-}\delta)
\rho = {1\over q}\rho (\varphi_{-}\delta)&\;\;(\varphi_{+}\delta) \rho =
{1\over q}\rho (\varphi_{+}\delta) - {(q^2-1)\over q}\varphi_{+}
(\rho\delta)\\ (\rho\delta) \varphi_{-} = \varphi_{-}
(\rho\delta)&\;\;(\varphi_{-}\delta) \varphi_{-} = \varphi_{-}
(\varphi_{-}\delta)&\;\;(\varphi_{+}\delta) \varphi_{-} = \varphi_{-}
(\varphi_{+}\delta)\\ (\rho\delta) \varphi_{+} = \varphi_{+}
(\rho\delta)&\;\;(\varphi_{-}\delta) \varphi_{+} = \varphi_{+}
(\varphi_{-}\delta)&\;\;(\varphi_{+}\delta) \varphi_{+} = \varphi_{+}
(\varphi_{+}\delta) \end{array} \eqno(27)
$$ are consistent with the C.R.
(3) and are simpler than the C.R. of the differential calculus (6).  The
 $U_{q}(SL(2,R))$ algebra is the condition of the compatibility of the
relations (25)
$$
\begin{array}{c}
X^A(q^2\delta_{R^1}\delta_{R^3} - q^{-2}\delta_{R^3}\delta_{R^1} =
(q^2+1)\delta_{R^3})\\
X^A(q^2\delta_{R^2}\delta_{R^1} - q^{-2}\delta_{R^1}\delta_{R^2} =
(q^2+1)\delta_{R^2})\\
X^A(\delta_{R^3}\delta_{R^2} - q^2\delta_{R^2}\delta_{R^3} =
\delta_{R^1}) \end{array} \eqno(28) $$

{\bf 5.Equations of motion}\\
  We use the extremum principle of the action to obtain the
equations of motion  and
we must to commute the variations of fields and their derivatives on the
right or on the left side. We can use both variation $dX^A$ and
$\delta X^A$ to do this. The C.R.  of the differential calculus on the
$SL_{q}(2,R)$ group are insufficient to do this. Therefore, we need in
the differential calculus on the Lagrangian manifolds $(\rho
,\varphi_{\pm},\dot\rho ,\dot{\varphi_{\pm}} ,\acute{\rho} ,
\acute{\varphi_{\pm}} )$.  This is not the quantum group manifold and we
can not use the formalism of 1-forms. We can
require, that the Lagrangian equation of motion be coincident with the
conservation law $\partial_{\mu}{\omega}^{\mu}=0$ for Lagrangian with
$SL_{q}(2,R)$ group symmetry. At last, we can investigate the 1-
dimensional $\sigma$- models.
 The variational calculus is more suitable to obtain the equations of
motion. The C.R. between the $X^A$, $\dot{X^A}$, $\acute{X^A}$ and
$R^n$,
$\dot{R^n}$, $\acute{R^n}$ can obtain by differentiating the relations
(24).  $$ \begin{array}{ccc} \dot{\rho}R^1=q^2 R^1
\dot{\rho}&\;\dot{\rho}\dot{R^1}=q^2\dot{R^1}\dot{\rho}&\;
\rho\dot{R^1}=\dot{R^1}\rho\\
\dot{\rho}R^2=qR^2\dot{\rho}&\;\dot{\rho}\dot{R^2}=q\dot{R^2}\dot{\rho}&\;
\rho\dot{R^2}=q\dot{R^2}\rho\\
\dot{\rho}R^3=qR^3\dot{\rho}&\;\dot{\rho}\dot{R^3}=q\dot{R^3}\dot{\rho}&\;
\rho\dot{R^3}=q\dot{R^3}\rho
\end{array}
\eqno(29)
$$
 The derivatives of $\varphi_{\pm}$ commute with the derivatives of
$R^n$.
\\ {\bf 5. One dimensional $\sigma$- model on
the quantum plane ($C_{q}(2|0)$.}\\ The differential
calculus on the $C_{q}(2|0)$ is coincide with the differential calculus
on the Borel subgroup of $SL_{q}(2,C)$ and can be obtained from the
differential calculus on the $SL_{q}(2,C)$ by surjection:
$\pi$:$SL_{q}(2,C)\rightarrow B_{L}$ such that $\pi(b)= 0$.
$$
g=\left( \begin{array}{cc} x&0\\ y&x^{-1} \end{array} \right),\;\;
\begin{array}{cc} {xy=qyx}&\dot yy=q^{-2}y\dot y \\ \dot xx=q^{-2}x\dot
x
&\dot xy=q^{-1}y\dot x \end{array} \eqno(30) $$ $$ \omega= \left(
\begin{array}{cc} x^{-1}dx&0\\ xdy-qydx&-q^{2}x^{-1}dx \end{array}
\right),\; \dot yx=q^{-1}x\dot y - {(q^{2}-1)\over q^{2}}y\dot x $$ In
term of the variables $\rho,\varphi_{\pm}$ $$g= \left( \begin{array}{cc}
\rho&0\\ \varphi_{+}\rho&\rho^{-1} \end{array} \right);\;\; \omega=
\left(
\begin{array}{cc}
\rho^{-1}d\rho&0\\
{1\over q}\rho^{2}d\varphi_{+}&-q^{2}\rho^{-1}d\rho
\end{array}
\right)
\eqno(31)
$$
$$
L={1\over 2}Tr_{q}(\omega_{\mu}\omega^{\mu})={q^{4}(q^{2}+1)\over 2}
\rho^{-2}{\dot\rho}^{2}
$$
 The equation of motion  ${\dot\omega}^{1}=
 \rho^{-1}\ddot{\rho}-q^{2}\rho^{-1}{\dot\rho}^{2}=0$ will
coincide with Lagrangian equation , if we impose the C.R.
$\delta\dot\rho\dot\rho={1\over
q^{2}}\dot\rho\delta\dot\rho$.
 The classical solution of this equation is
$$\rho=\alpha\exp(\beta t),\;\;\
\alpha\beta=q^{2}\beta\alpha
\eqno(32)
$$ and C.R.
$$\rho(t)\rho(t^{\prime})=\rho(q^{2}t^{\prime})\rho({1\over
q^{2}}t),\;\;
\rho(t)\rho(t^{\prime})= \exp[q^{2}(q^{2}-1)\beta
(t-t^{\prime})]\rho(t^{\prime})\rho(t)
\eqno(33)
$$. There are
$4\times 4$ matrix representations of $\alpha, \beta$ such,that
$det_{q}\alpha= 0$ or $det_{q}\beta= 0$.Therefore, we can
rewrite this Lagrangian as a $4\times 4$ matrix model for the commuting
fields.In conclusion, note that 2-dimensional $\sigma$-model on
the guantum plane $$L={q^{4}(q^{2}+1)\over
2}\rho^{-2}\partial_{\mu}\rho\partial^{\mu}\rho
\eqno(34)
$$ leads to the C.R.
$\delta\acute\rho\acute\rho={1\over
q^{2}}\acute\rho\delta\acute\rho$ and the equation of motion
$\partial_{\mu}\partial^{\mu}\rho-q^{2}\rho^{-1}\partial_{\mu}
\rho\partial^{\mu}\rho=0,\;\;\mu=1,2$.

\vspace{0.20cm}

  I would like to thank J.Wess, V.Dobrev, J,Lukierski, A.Isaev,
P.Pyatov,
B.Zupnik,\\ V.Lyakhovsky, A.Akulov for stimulating discussions.

  This work was supported in part by the Ukrainian Ministry of Science
and
technology, grant  2.5.1/54, by grants INTAS 93-633 (Extension) and
INTAS 93-127 (Extension).
\begin{center}
{\bf References}
\end{center}
\begin{enumerate}
\item L.Faddeev, N.Reshetikhin, L.Takhtajan, {\sl Algebra i
Analiz} { \bf 1} (1989), p.178
\item Akulov V.P., Gershun V.D., Gumenchuk A.I., {\sl JETP
Lett.}, {\bf 56} (1992), p.177
\item Akulov V.P., Gershun V.D., Gumenchuk A.I., {\sl JETP
Lett.}, {\bf 58} (1993), p.474
\item V.Akulov, V.Gershun, {\sl q-alg 9509030}
\item Akulov V.P., Gershun V.D.,{\sl Proc. Int. conf. on math, phys.
(Rahov, Ukraine, 1995)}, (in press)
\item L.D.Fadeev, {\sl Commun. Math, Phys.}, {\bf 132} (1990), p.131
\item A.Alekseev and S.Shatashvili, {\sl Commun. Math. Phys.}, {\bf 133}
(1990), p.353
\item Arefeva I.Y. and Volovich I.V., {\sl Phys. Lett.}, {\bf 264B}
(1991), p.62
\item Y.Frishman, J.Lukierski, W.J.Zakrzewski, {\sl J.Phys.A:Math.Gen},
{\bf 26} (1993), p.301
\item V.D.Gershun, {\sl Proc. Xth Int. conf. "Problems of Quantum Field
Theory" (Alushta, Ukraine,1996)}, JINR, Dubna, (1996), p.119
\end{enumerate}

\end{document}